\documentclass[a4paper,12pt,draft]{amsart}
\usepackage{amsmath, amssymb,braket}
\usepackage{amsthm}
\usepackage[vcentermath,enableskew]{youngtab}
\usepackage{bm}
\newcommand{\partitionof}{\vdash}

\newcommand{\YYoungsub}[1]{\boldsymbol{S}_{#1}}
\newcommand{\colsub}[2]{C_{#1,#2}}
\newcommand{\CColsub}[3]{\boldsymbol{C}_{#1,#2,#3}}
\newcommand{\ZZsub}[3]{\boldsymbol{Z}_{#1,#2,#3}}
\newcommand{\SSymmetrizer}[1]{\boldsymbol{A}_{#1}}
\newcommand{\stdtab}[1]{\mathcal{T}_{#1}}
\newcommand{\indextab}[1]{\mathcal{I}_{#1}}
\newcommand{\stdtabp}[1]{\mathcal{T}'_{#1}}
\newcommand{\stdtabpp}[1]{\mathcal{T}''_{#1}}
\newcommand{\indextabpair}[1]{\mathcal{J}_{#1}}

\newcommand{\xstdtab}[1]{\mathcal{T}_{#1}}
\newcommand{\xindextab}[1]{\mathcal{I}_{#1}}
\newcommand{\xstdtabp}[1]{\mathcal{T}'_{#1}}
\newcommand{\xindextabpair}[1]{\mathcal{J}_{#1}}

\newcommand{\ZZ}{\mathbb{Z}}

\newcommand{\NN}{\mathbb{N}}
\newcommand{\CC}{\mathbb{C}}

\newcommand{\numof}[1]{\left| #1 \right|}

\newcommand{\induce}[3]{\operatorname{Ind}_{#1}^{#2}{#3}}
\newcommand{\Ann}{\operatorname{Ann}}

\newcommand{\disjointunion}{\amalg}
\newcommand{\defit}[1]{{\em #1}}

\theoremstyle{plain}
\newtheorem{thm}{Theorem}[section]
\newtheorem{lemma}[thm]{Lemma}

\newtheorem{theorem}[thm]{Theorem}

\theoremstyle{remark}
\newtheorem{remark}[thm]{Remark}

\newtheorem{example}[thm]{Example}

\theoremstyle{definition}
\newtheorem{definition}[thm]{Definition}

\newtheorem{notation}[thm]{Notation}

\allowdisplaybreaks[3]

\begin{document}

\title[An extended Schur's lemma]{
An extended Schur's lemma and its application
}
\author[Numata, y.]{Numata, yasuhide}
\address{Hokkaido university}
\email{nu@math.sci.hokudai.ac.jp}
%\urladdr{http://www.math.sci.hokudai.ac.jp/\textasciitilde nu/}

\begin{abstract}
The Springer modules have a combinatorial property called 
``coincidence of dimensions,''
i.e., the Springer modules are naturally decomposed into submodules with common dimensions.
Morita and Nakajima proved the property by giving modules with
 common dimensions whose induced modules are isomorphic to  the
 submodules of Springer modules.
They proved that the induced modules are isomorphic to the submodules, by
showing the coincidence of their characters.
Our aim is to construct isomorphisms between the induced modules
 and their corresponding submodules in a combinatorial manner.
For this purpose, we show lemmas, which are equivalent to the classical
Schur's lemma in special cases.
We also give  a procedure to construct isomorphisms, and
explicitly construct isomorphisms in the case of the Springer modules
corresponding to Young diagrams of two rows.
\end{abstract}

\maketitle

\section{Introduction}
Let $G$ be a finite group.
In some $\ZZ$-graded $G$-modules $R=\bigoplus_{d} R^{d}$, we have 
a phenomenon called ``coincidence of dimensions,'' i.e.,
some integers $l$  satisfy
the equations
\begin{gather*}
\dim \bigoplus_{i\in\ZZ} R^{i l+k} = \dim \bigoplus_{i\in\ZZ} R^{i l+k'}
\end{gather*}
for all $k$ and $k'$.
We call a datum $(H(l),\Set{Z(k;l)})$ 
a \defit{representation-theoretical presentation}  for the phenomenon
if a subgroup $H(l)$ of $G$ and $H(l)$-modules $Z(k;l)$ satisfy
\begin{align*}
\bigoplus_{i\in\ZZ} R^{i l+k} \simeq \induce{H(l)}{G}{Z(k;l)},  && \dim Z(k;l)=\dim Z(k';l)
\end{align*}
for all $k$ and $k'$, where $\induce{H(l)}{G}{Z(k;l)}$ denotes the
induced module.

Morita and Nakajima 
gave  representation-theoretical presentations to
coincidences of dimensions of the Springer modules
$R_\mu$ \cite{MN,m,mn,Mo}. 
The Springer modules $R_\mu$ are graded algebras parametrized by
 partitions $\mu \partitionof m$.
As $S_m$-modules, $R_\mu$ are  
isomorphic to the cohomology rings
of the variety of flags
 fixed by a unipotent matrix 
with Jordan blocks of type $\mu$.
  (See \cite{HotSpring, Sp76, Sp78}. 
See also \cite{dp, tani} for algebraic construction.)
We recall the case where $\mu$ is an
$l$-partition, 
where an $l$-partition means a partition whose
multiplicities are divisible by $l$.
Let $R_\mu(k;l)$ denote
the submodule $\bigoplus_{i\in\ZZ} R_\mu^{i l+k}$ of 
the Springer module $R_\mu$.
In this case, we have 
$\dim R_\mu(k;l)=\dim R_\mu(k';l)$ for all $k$ and $k'$, i.e.,
$R_\mu$ has a coincidence of dimensions.
Let $H_{\mu}(l)$ be the semi-direct product 
$S_\mu \rtimes \colsub{\mu}{l}$ of 
the Young subgroup $S_\mu$ and an $l$-th cyclic group
$\colsub{\mu}{l}=\Braket{a_{\mu,l}}$. 
(See Section \ref{notationsec} for their definitions.)
For $k\in\ZZ$, let $ Z_\mu(k;l): H_{\mu}(l) \to \CC^{\times} $ 
denote one-dimensional representations  of $H_{\mu}(l)$
mapping $a_{\mu,l}$ to $\zeta_l^k$ and $\sigma\in S_m$ to $1$,
where $\zeta_l$ denotes a primitive $l$-th root of unity.
Then  $(H_\mu (l),\Set{Z_\mu(k;l)})$ 
is a representation-theoretical presentation for this case, i.e.,
 $R_{\mu}(k;l) \simeq \induce{H_\mu(l)}{S_m}{Z_\mu(k;l)}$ for all $k$.
To prove it,
Morita and Nakajima described the values of the Green polynomials 
at roots of unity, and showed that the characters of the submodules
 $ R_{\mu}(k;l)$
coincide with those of the induced modules
$\induce{H_\mu(l)}{S_m}{Z_\mu(k;l)}$ in \cite{mn}.

Our motivation in this paper is to construct isomorphisms 
from the induced module 
$\induce{H_\mu(l)}{S_m}{Z_\mu(k;l)}$ 
to the submodule $R_{\mu}(k;l)$ in a combinatorial manner.

For this purpose, 
we consider the problem in a more general situation.
In Section \ref{mainsec}, 
we consider  
the induced module  $\induce{H}{G}{Z}$
of a one-dimensional $H$-module $Z$
and  another realization $R$ of  $\induce{H}{G}{Z}$,
where  $H$ is a subgroup of a finite group $G$.
We show lemmas about
the condition for an element of $R$
to generate $R$.
In special cases,
these lemmas are  equivalent to the classical Schur's lemma.
We give  a necessary and sufficient  condition 
for a map 
from  $\induce{H}{G}{Z}$ to  $R$ 
to be an isomorphism.
We also give a procedure to construct isomorphisms.
In Subsection \ref{appsec}, we apply 
the main results to the case of representation-theoretical presentations
for the Springer modules corresponding to $l$-partitions $\mu$.
In Subsection \ref{tworowsec}, we explicitly give isomorphisms in the
case of $\mu=(n,n)$

\section{Main Results}\label{mainsec}
In this section, we consider a subgroup $H$ of a finite group  $G$
and a one-dimensional representation $\zeta$ of  $H$.
Let $M$ be the induced module of $\zeta$, and $R$ a $G$-module
isomorphic to $M$.
We show a condition 
for a map from $M$ to $R$ to be an isomorphism of $G$-modules.
We prove Theorems \ref{modschulemmaIII}, \ref{mainlemma} and
\ref{lemmagroupring}  in Subsection \ref{proofsec}.

\begin{notation}
Throughout this section, we use the following notation and assumption.
Let $K$ be a field. We consider representations over $K$.
We assume the complete reducibility of representations.
Let $G$ be a finite group,
    $\varepsilon$  the unit of $G$,
    $H$    a subgroup such that $\numof{H}\in K^{\times}$,
    $Z$    a one-dimensional vector space over $K$, 
    $\zeta:H\to K^{\times}$ a representation acting on $Z$,
and
    $\induce{H}{G}{Z}$  the induced module of $Z$.
\end{notation}

\begin{definition}\label{genericdef} 
For a $G$-module $R$,
 $f \in R$ is said to be \defit{generic in $R$}
if there exists an irreducible decomposition
$R=\bigoplus_{\lambda} \bigoplus_{P \in \indextab{\lambda}} S^{P}$
satisfying the following:
\begin{itemize}
\item $S^P$ is isomorphic to $S^{P'}$ if and only if  $P, P'\in      \indextab{\lambda}$ 
for some $\lambda$.
\item 
      For $P\in \indextab{\lambda}$, let $f^{(P)}$ denote the image of
      the projection of $f$ to the irreducible component $S^{P}$.
      For $P, P'\in \indextab{\lambda}$, let $\psi_{P,P'}$ be
      an isomorphism   from $S^{P}$ to $S^{P'}$.  
       Then $\Set{ \psi_{P,P'} f^{(P)}| P \in \indextab{\lambda} } \subset S^{P'}$
is $K$-linearly independent. 
\end{itemize}
\end{definition}

\begin{remark}
By Schur's lemma, Definition \ref{genericdef} does not depend on
choices of $\psi_{P,P'}$.
\end{remark}

We can describe a  necessary and sufficient condition for an element of
$R$ to generate $R$ as a $G$-module.

\begin{theorem}[An extended Schur's lemma]\label{modschulemmaIII}
For a $G$-module $R$,
$F \in R$ generates $R$ if and only if $F$ is generic in $R$.
\end{theorem}
\begin{remark}
In the case where $R = S \oplus S$ and $S$ is a simple $G$-module,
this theorem is equivalent to the classical Schur's lemma.
\end{remark}

\begin{definition} 
Let $M$ and $R$ be $G$-modules.
Let $e \in M$ generate $M$ as a $G$-module.
For $f \in R$, 
we define $\varphi_{e,f}$ to be the $K[G]$-linear map such that
\begin{gather*}
\varphi_{e,f}: M \ni e \mapsto  f \in R
\end{gather*}
whenever it is well-defined,
where $K[G]$ denotes the group ring of $G$.
\end{definition}

\begin{definition} 
For a representation $\zeta: H \to K^{\times}$,
we define an element $\SSymmetrizer{\zeta}$ 
of the group ring $ K [H]$ to be 
$\frac{1}{\numof{H}}\sum_{\sigma\in H} \zeta(\sigma^{-1}) \sigma $.
\end{definition}

We can describe  a necessary and sufficient  condition 
for a map from  $\induce{H}{G}{Z}$ to  $R$ 
to be an isomorphism.

\begin{theorem}\label{mainlemma}
Let $R$ be a realization of $\induce{H}{G}{Z}$,
 $e\in Z$ a nonzero element, 
and
$\phi$  a map from $\induce{H}{G}{Z}$ to $R$.
Then $\phi$ is an isomorphism from $\induce{H}{G}{Z}$ to $R$
if and only if
$\phi=\varphi_{\varepsilon \otimes e, \SSymmetrizer{\zeta} f}$
for some $\SSymmetrizer{\zeta} f$ which is generic in $R$.
\end{theorem}

\begin{theorem}\label{lemmagroupring}
The induced module $\induce{H}{G}{Z}$ 
is isomorphic to $K[G]\SSymmetrizer{\zeta}$ as $G$-modules.
For a nonzero element $e\in Z$,
the map $\varphi_{\varepsilon \otimes e, \SSymmetrizer{\zeta}}$  is an isomorphism from $\induce{H}{G}{Z}$ to $K[G]\SSymmetrizer{\zeta}$.
\end{theorem}

By Theorem \ref{mainlemma}, constructing isomorphisms is equivalent to
finding  elements $f$ such that $\SSymmetrizer{\zeta} f$
are generic.
By the definition of generic elements, we have 
the following procedure to construct such elements.
\begin{theorem}\label{procelemma1} 
For all elements $f$ given by the following procedure,
$\SSymmetrizer{\zeta} f$
are generic in $R$ $:$
\begin{enumerate}
\item Fix an irreducible decomposition 
       $\bigoplus_{\lambda} \bigoplus_{P \in \xindextab{\lambda}} S^{P}$
      of $R$
      such that $S^{P}$ is isomorphic to $S^{P'}$
      if and only if  $P, P' \in \xindextab{\lambda}$ for some $\lambda$. 
\item Fix  an element $P_{\lambda}$  of $\xindextab{\lambda}$ for each
      $\lambda$. 
      Let $S^{\lambda}=S^{P_\lambda}$. 
\item Fix  an isomorphism $\psi_P : S^\lambda \to S^P$ for 
      each $P \in \xindextab{\lambda}$. 
\item Fix  a basis $\Set{\Delta_{Q}|Q\in\xstdtab{\lambda}}$ of $S^{\lambda}$. 
      Let $\Delta_Q^P$ be the image $\psi_P(\Delta_Q)$.
\item Fix  a maximal subset $\xstdtabp{\lambda}$  of  $\xstdtab{\lambda}$
      such that 
      $\Set{ \SSymmetrizer{\zeta} \Delta_{Q} | Q \in \xstdtabp{\lambda} }$
      is $K$-linearly independent. 
\item Fix a total ordering on $\xstdtabp{\lambda}$.
      Let $\Set{Q_1<Q_2<\cdots} = \xstdtabp{\lambda} $.
\item Fix a total ordering on $\xindextab{\lambda}$.
      Let $\Set{P_1<P_2<\cdots<P_d} = \xindextab{\lambda} $. 
\item Let $\xindextabpair{\lambda} = \Set{(P_i,Q_i)| i=1, \ldots , d}$.
\item Let $f$ be 
       $\sum_{\lambda}\sum_{(P,Q)\in \xindextabpair{\lambda}} \Delta^{P}_{Q}$.
\end{enumerate}
\end{theorem}

\begin{remark}
By definition,
any element $f\in R$ such that $\SSymmetrizer{\zeta} f$ is generic in $R$
is obtained from the procedure in Theorem \ref{procelemma1}.
\end{remark}

\subsection{Proof of Main Results}\label{proofsec}
In this section, we prove the main results.
First, in \ref{necessarysec},  
we prove Theorem \ref{modschulemmaIII}, an extended Schur's lemma.
Next, in \ref{welldefsec}, we show
a necessary  and sufficient condition
for  $\varphi_{e,f}$ to be injective and well-defined as
 a homomorphism of $G$-modules.
Finally we prove  Theorems \ref{mainlemma} and \ref{lemmagroupring}.

\subsubsection{Proof of Theorem \ref{modschulemmaIII}}\label{necessarysec}
Here we prove 
Theorem \ref{modschulemmaIII}, 
which gives a necessary and sufficient condition for the surjectivity of $\varphi$.

\begin{lemma}\label{modschulemmaI}
For a simple $G$-module $S$ and a subset $\Set{f_1,\ldots,f_n} \subset S$,
\begin{gather*}
  f_1\oplus \cdots \oplus f_n \in \bigoplus_{i=1}^{n} S
\end{gather*}
 generates $\bigoplus_{i=1}^{n} S$
if and only if $\Set{f_1,\ldots,f_n}$ is $K$-linearly independent.
\end{lemma}
\begin{proof}
First we prove the case where $n=2$.

The case where $f_1=0$ or $f_2=0$ is clear.

Suppose that $f_1\neq 0$ and $f_2\neq 0$.
There exists $P_1\in K[G]$ such that 
$P_1 (f_1 \oplus f_2) = f_1 \oplus 0$
in $S\oplus S$
if and only if $K[G] \Set{f_1\oplus f_2} = S \oplus S$.

If $f_1 = \alpha f_2$, then $\Ann(f_1) = \Ann(\alpha f_2) = \Ann(f_2)$.
Hence there does not exist such $P_1$. 

Conversely we consider the case where 
 there does not exist such $P_1$. 
Since this implies $P_1 f_1=0$ for all $P_1 \in \Ann(f_2)$,
we have $\Ann(f_2) \subset \Ann(f_1)$.
Hence $\varphi_{f_2,f_1}$ is a well-defined endomorphism of $S$.
By the classical Schur's lemma, $f_1$ and $f_2$ are $K$-linearly dependent.

Thus we have the lemma for $n=2$.

Next we prove the lemma in the case where $n=k$.

Assume that $f=f_1 \oplus \cdots \oplus f_{k-1} \in \bigoplus_{i=1}^{k-1} S$ generates
 $\bigoplus_{i=1}^{k-1} S$.
In this case, there exists $P_i$ such that
\begin{gather*}
P_i f =0\oplus \cdots \oplus 0 \oplus f_i \oplus 0 \oplus \cdots \oplus
 0
\end{gather*}
 for each $i$.

If  $\Set{f_1,\ldots ,f_{k-1} ,f_k}$ is  $K$-linearly dependent, 
then 
\begin{gather*}
\bigcap_{i=1}^{k-1} \Ann (f_i) \subset \Ann(f_k).
\end{gather*}
Hence $K[G]\Set{f \oplus f_k}\simeq K[G]\Set{f} = \bigoplus_{i=1}^{k-1} S$.

Conversely,
we consider the case where $K[G]\Set{f \oplus f_k} \simeq \bigoplus_{i=1}^{k-1} S$.
In this case, $\Ann(g)\subset\bigcap_{i=1}^{k-1} \Ann(f_i)$.
Since 
$P_i (f \oplus f_k) = 0\oplus \cdots \oplus 0 \oplus f_i \oplus 0 \oplus \cdots \oplus 0 \oplus P_i f_k,  $
\begin{gather*}
K[G]\Set{P_i (f \oplus f_k)} \simeq K[G] \Set{ f_i\oplus P_i f_k}\simeq S.  
\end{gather*}
Since it follows from the case where $n=2$ that
there exists  $\alpha_i \in K$ such that 
$P_i f_k = \alpha_i f_i$, 
\begin{gather*}
\sum_{i=1}^{k-1} P_i f_k = \sum_{i=1}^{k-1} \alpha_i f_i.
\end{gather*}
Since $(1- \sum_{i=1}^{n} P_i)(f \oplus f_k) = 0 \oplus \cdots \oplus 0 \oplus  (f_k - \sum_{i=1}^{n} \alpha_i f_i)$,
\begin{gather*}
f_k-\sum_{i=1}^{k-1} \alpha_i f_i = 0.
\end{gather*}
Hence $\Set{f_1,\ldots ,f_{k-1} ,f_k}$ is  $K$-linearly dependent.

Thus we have the lemma in the case where $n=k$.
\end{proof}

\begin{remark}
In the case where $n=2$,
Lemma \ref{modschulemmaI} is equivalent to 
the classical Schur's lemma.
\end{remark}

\begin{lemma}\label{modschulemmaII}
Let  $\Set{S^{\lambda}}$ be a family of  simple $G$-modules
which are not isomorphic to one another.
For subsets 
$\Set{f^{\lambda}_1,\ldots,f^{\lambda}_{m_\lambda}} \subset S^{\lambda}$,
\begin{gather*}
\bigoplus_{\lambda} f^{\lambda}_1 \oplus \cdots \oplus f^{\lambda}_{m_\lambda}
\in \bigoplus_{\lambda} \bigoplus_{i=1}^{m_\lambda} S^{\lambda}
\end{gather*}
 generates $\bigoplus_{\lambda}\bigoplus_{i=1}^{m_\lambda} S^{\lambda}$
if and only if each subset $\Set{f^{\lambda}_1,\ldots,f^{\lambda}_{m_\lambda}}$ is $K$-linearly independent.
\end{lemma}
\begin{proof}
By Lemma \ref{modschulemmaI},
$f^{\lambda}_1 \oplus \cdots \oplus f^{\lambda}_{m_\lambda}$
generates $\bigoplus_{i=1}^{m_\lambda} S^{\lambda}$
if and only if
$\Set{f^{\lambda}_1,\ldots,f^{\lambda}_{m_\lambda}}\subset S^{\lambda}$
 is $K$-linearly independent.
Since $S^{\lambda}$ is not isomorphic to $S^{\lambda'}$ 
for $\lambda\neq\lambda'$,
\begin{gather*}
\bigoplus_{\lambda} f^{\lambda}_1 \oplus \cdots \oplus f^{\lambda}_{m_\lambda}
\in \bigoplus_{\lambda} \bigoplus_{i=1}^{m_\lambda} S^{\lambda}
\end{gather*}
 generates $\bigoplus_{\lambda}\bigoplus_{i=1}^{m_\lambda} S^{\lambda}$
if and only if each subset $\Set{f^{\lambda}_1,\ldots,f^{\lambda}_{n_\lambda}}\subset S^{\lambda}$ is $K$-linearly independent.
\end{proof}

\begin{proof}[Proof of Theorem $\ref{modschulemmaIII}$]
This follows from the definition of generic elements and Lemma \ref{modschulemmaII}.
\end{proof}

\subsubsection{Proofs of Theorems \ref{mainlemma} and \ref{lemmagroupring}}
\label{welldefsec}
Here we show Lemma \ref{symmetrizerlemma}, which 
gives 
a necessary  and sufficient  condition
for $\varphi_{e,f}$ to be injective
and well-defined as a homomorphism of $G$-modules.
Then we prove  Theorems \ref{mainlemma} and \ref{lemmagroupring}.

\begin{lemma}\label{symmetrizerlemma}
Let $G$ be a group, $H$ a subgroup of $G$, 
and $\zeta: G \to K^{\times}$ a one-dimensional
 representation.
For a $G$-module $V$,
\begin{align*}
\Set{ v \in V | \sigma v = \zeta(\sigma) v \text{ for all } \sigma \in H}
=\Set{ \SSymmetrizer{\zeta} v | v \in V }.
\end{align*}
\end{lemma}

\begin{proof}
For $v\in \Set{ v \in V | \sigma v = \zeta(\sigma) v \text{ for all } \sigma \in H}$,
\begin{align*}
\SSymmetrizer{\zeta} v 
&= \frac{1}{\numof{H}}\sum_{\sigma\in H} \zeta(\sigma^{-1}) \sigma v\\
& = \frac{1}{\numof{H}}\sum_{\sigma\in H} \zeta(\sigma^{-1}) \zeta(\sigma) v\\
& = \frac{1}{\numof{H}}\sum_{\sigma\in H} v\\
& = v.
\end{align*}
Hence
 $v\in \Set{ \SSymmetrizer{\zeta} v | v \in V }$.

For $\tau \in H$,
\begin{align*}
\tau \SSymmetrizer{\zeta} v &=\tau \frac{1}{\numof{H}}\sum_{\sigma\in H} \zeta(\sigma^{-1}) \sigma v\\
& = \frac{1}{\numof{H}}\sum_{\sigma\in H} \zeta(\sigma^{-1}) \tau \sigma v\\
& = \frac{1}{\numof{H}}\sum_{\sigma\in H}
 \zeta(\sigma^{-1})\zeta(\tau^{-1}) \zeta(\tau)\tau \sigma v \\
& = \frac{1}{\numof{H}}\sum_{\sigma'\in H} \zeta(\sigma'^{-1})
 \zeta(\tau)\sigma' v \\
& =  \zeta(\tau)\frac{1}{\numof{H}}\sum_{\sigma'\in H} \zeta(\sigma'^{-1})
\sigma' v \\
& = \zeta(\tau)\SSymmetrizer{\zeta} v. 
\end{align*}
Hence $\SSymmetrizer{\zeta} v \in \Set{ v \in V | \sigma v = \zeta(\sigma) v \text{ for all } \sigma \in H}$.
Thus we have
\begin{gather*}
\Set{ v \in V | \sigma v = \zeta(\sigma) v \text{ for all } \sigma \in H}
=\Set{ \SSymmetrizer{\zeta} v | v \in V }.
\end{gather*}
\end{proof}

\begin{remark}
For any $\zeta$, the operator $\SSymmetrizer{\zeta}$ is a projection,
 i.e.,
 $\SSymmetrizer{\zeta}^2 = \SSymmetrizer{\zeta}$.
\end{remark}

Now we prove Theorems $\ref{mainlemma}$ and  $\ref{lemmagroupring}$.

\begin{proof}[Proof of Theorem $\ref{mainlemma}$]
Let $\phi$ be a homomorphism from $\induce{H}{G}{Z}$ to $R$.
Since $Z$ is one-dimensional, $\phi$ is determined by 
$\phi(\varepsilon\otimes e)$.
By Lemma $\ref{symmetrizerlemma}$, 
the homomorphism
$\varphi_{\varepsilon\otimes e,\SSymmetrizer{\zeta} f}$  is well-defined.
Conversely, if we take $g$ which can not be written as $\SSymmetrizer{\zeta} f$,
then $\varphi_{\varepsilon\otimes e,g}$ is not an isomorphism.
Hence we have $\phi=\varphi_{\varepsilon\otimes e,\SSymmetrizer{\zeta} f}$ 
for some $f$.

Now $\varphi_{\varepsilon\otimes e,\SSymmetrizer{\zeta} f}$ is
 well-defined, and 
$\induce{H}{G}{Z}$ is isomorphic to $R$.
Since $R$ is finite dimensional, 
it is enough to prove the surjectivity 
$\varphi_{\varepsilon\otimes e,\SSymmetrizer{\zeta} f}$.
By Theorem \ref{modschulemmaIII},
 $K[G]\Set{\SSymmetrizer{\zeta} f}=R$ if and only if $\SSymmetrizer{\zeta} f$ is generic in $R$.
Hence 
$\varphi_{\varepsilon\otimes e,\SSymmetrizer{\zeta} f}$ is surjective
 if and only if $\SSymmetrizer{\zeta} f$ is generic in $R$.
Therefore we have  Theorem \ref{mainlemma}.
\end{proof}

\begin{proof}[Proof of Theorem $\ref{lemmagroupring}$]
By Lemma \ref{symmetrizerlemma},
$K\SSymmetrizer{\zeta}$ is isomorphic to $Z$ as $H$-modules.
By definition, $\tau K\SSymmetrizer{\zeta}=\tau'K\SSymmetrizer{\zeta}$
if and only if $\tau H= \tau' H$.
Hence we have Theorem \ref{lemmagroupring}.
\end{proof}

\section{Application}
In this section, we apply the main results to the case of the Springer modules.

In Subsection \ref{notationsec}, we prepare some basic notation. 
In Subsection \ref{appsec}, we rewrite the main results for the case of
representation-theoretical presentations for
the Springer modules $R_\mu$  corresponding to 
$l$-partitions $\mu$.
Moreover, in Subsection \ref{tworowsec},
we explicitly describe isomorphisms in the case where $\mu=(n,n)$.

\subsection{Notation}\label{notationsec}
We identify a partition $\mu=(\mu_1\geq\mu_2\geq\cdots)$ of $m$ with its
Young diagram  $\Set{(i,j) \in \NN^2 |1 \leq j \leq \mu_i}$ with $m$ boxes.
If $\mu$ is a Young diagram with $m$ boxes, 
we write $\mu \partitionof m$ and identify a Young diagram $\mu$ with 
the array of $m$ boxes having left-justified rows with the $i$-th row
containing  $\mu_i$ boxes; for example,
\begin{gather*}
(2,2,1) = \yng(2,2,1)  \partitionof 5.
\end{gather*}
For an integer $l$, a Young diagram $\mu$ is called an \defit{$l$-partition}
if multiplicities $m_i=\numof{\Set{k|\mu_k=i}}$ of $i$
are divisible by $l$ for all $i$.

We call a map  $T$ a \defit{numbering on a Young diagram $\mu \partitionof m$} 
if
$T$ is a bijection $\mu \ni (i,j) \mapsto T_{i,j} \in \Set{1,\ldots , m}$.
We call a map  $T$ a \defit{semi-standard tableau on a Young diagram 
$\mu\partitionof m$ with weight $w=(w_1,w_2,\ldots)$} 
if
$T:\mu \ni (i,j) \mapsto T_{i,j} \in \Set{1,\ldots , m}$
satisfies $T_{i,j}<T_{i+1,j}$, $T_{i,j}\leq T_{i,j+1}$ and
$\numof{T^{-1}(\Set{k})}=w_k$ for all $i,j,k$.
A numbering which is a semi-standard tableau is called a \defit{standard
tableau}.
We identify a map $T$ from a Young diagram $\mu$ to $\NN$ 
with a diagram putting each $T_{i,j}$ on
the box in the $(i,j)$ position; for example, 
\begin{gather*}
\young(23,41,5).
\end{gather*}
For $\mu \partitionof m$,
 $t_\mu$ denotes the numbering that maps
$(t_\mu)_{i,j}=j+\sum_{k=0}^{i-1}\mu_k$, i.e.,
a numbering obtained by putting numbers from $1$ to $m$ on boxes of $\mu$
from left to right in each row, starting in the top row and moving
to the bottom row.
For example, 
\begin{gather*}
t_{\tiny\yng(2,2,1)}=\young(12,34,5) \ .
\end{gather*}

For $\mu \partitionof m$, $S_\mu$ denotes the Young subgroup 
\begin{gather*}
S_{\Set{1,2,\ldots,\mu_1}} \times
 S_{\Set{\mu_1 + 1 , \mu_1 + 2,\ldots,\mu_1 + \mu_2 }} \times 
\cdots .
\end{gather*}
We define the element $\YYoungsub{\mu}$ of the group ring $\CC[S_m]$
to be 
\begin{gather*}
\frac{1}{\numof{S_\mu}} \sum_{\sigma\in S_\mu} \sigma .
\end{gather*}

Let $T$ be a numbering on an $l$-partition  $\mu \partitionof m$.
We define $a_{T,l}$ to be the product 
\begin{gather*}
\prod_{(l i+1,j)\in\mu}(T_{l i+1,j}, T_{l i+2,j},\ldots, T_{l i+l,j})
\end{gather*}
of $m/l$ cyclic permutations of
length $l$.
For example, $a_{t_{(2,2,1,1)},2}=(13)(24)(56)$.
We write $a_{\mu,l}$ for $a_{t_{\mu},l}$.
We define $\colsub{\mu}{l}$ to be the $l$-th cyclic group $\Braket{a_{\mu,l}}$.
We define  $H_\mu(l)$ to be the semi-direct product
$S_\mu \rtimes \colsub{\mu}{l}$.
For $k\in\ZZ$, let $ Z_\mu(k;l): H_{\mu}(l) \to \CC^{\times} $ 
denote one-dimensional representations  of $H_{\mu}(l)$
that maps $a_{\mu,l}$ to $\zeta_l^k$ and $\sigma\in S_m$ to $1$,
where $\zeta_l$ denotes a primitive root of unity.

We  define the element $\CColsub{\mu}{k}{l}\in\CC[S_m]$ to be
\begin{gather*}
\frac{1}{l}\sum_{j \in \ZZ/l\ZZ} \zeta_l^{-k j} a_{\mu,l}^{j}. 
\end{gather*}
We also define $\ZZsub{\mu}{k}{l} \in \CC[S_m]$ to be
\begin{gather*}
\CColsub{\mu}{k}{l}\YYoungsub{\mu}.
\end{gather*}

\begin{remark}
By definition, $\SSymmetrizer{1_{S_\mu}} = \YYoungsub{\mu}$,
and $\SSymmetrizer{Z_{\mu}(k,l)} = \ZZsub{\mu}{k}{l}$,
where $1_{S_\mu}$ denotes the trivial representation of $S_\mu$.
\end{remark}

For $\mu \partitionof m$,
we define an $S_m$-module $R_{\mu}(k;l)$ to be
\begin{gather*}
\bigoplus_{i\in\ZZ} R_{\mu}^{i l+k},
\end{gather*}
where $R_{\mu}^{i l+k}$ is the component of degree $(i l + k)$ of the Springer
 module $R_\mu$.

\begin{remark}
By \cite{MN}, 
$\induce{H_\mu(l)}{S_m}{Z_\mu(k;l)}$, $\CC[S_m]\ZZsub{\mu}{k}{l}$ and $R_{\mu}(k;l)$
 are isomorphic to one
another.
\end{remark}

\subsection{Application}\label{appsec}
In this section, 
we apply the main results to the case of 
the Springer modules $R_\mu$
corresponding to $l$-partitions $\mu$.
We rewrite the main results for $\induce{H_\mu(l)}{S_m}{Z_\mu(k;l)}$, $\CC[S_m]\ZZsub{\mu}{k}{l}$, and $R_{\mu}(k;l)$.

\begin{theorem}\label{isothm1}
Let $\mu$ be an $l$-partition  of $m$.
Then 
$\varphi_{\ZZsub{\mu}{k}{l}, \ZZsub{\mu}{k}{l} f}$ is 
an isomorphism from $\CC[S_m]\ZZsub{\mu}{k}{l}$ to $R_\mu(k;l)$
if and only if $\ZZsub{\mu}{k}{l} f$ is generic in $R_\mu(k;l)$.
\end{theorem}

\begin{theorem}\label{isothm2}
Let $\mu$ be an $l$-partition of $m$.
Then  $\varphi_{\varepsilon\otimes e,  \ZZsub{\mu}{k}{l}}$ is 
an isomorphism from $\induce{H_\mu(l)}{S_m}{Z_\mu(k;l)}$ to $\CC[S_m]\ZZsub{\mu}{k}{l}$.
\end{theorem}

\begin{theorem}\label{isothm4}
Let $\mu$ be an $l$-partition of $m$.
Then  $\varphi_{\varepsilon\otimes e , \ZZsub{\mu}{k}{l} f}$ is 
an isomorphism from $\induce{H_\mu(l)}{S_m}{Z_\mu(k;l)}$ to $R_\mu(k;l)$
if and only if $\ZZsub{\mu}{k}{l}  f$ is generic in $R_\mu(k;l)$.
\end{theorem}

\begin{remark}
It also follows from  the proofs of the main results
that any isomorphism between these modules is of  such a form.
\end{remark}

We also have a procedure to get generic elements.

\begin{theorem}\label{procelemma} 
For all elements $f$ given by the following procedure,
$\ZZsub{\mu}{k}{l} f$ are generic in $R_{\mu}(k;l)$ $:$
\begin{enumerate}
\item Let 
$\bigoplus_{\lambda \partitionof m} \bigoplus_{P \in \indextab{\lambda}} S^{P}$
      be an irreducible decomposition of $R_{\mu}(k;l)$
such that $S^{P}$ is isomorphic to the Specht module $S^{\lambda}$
for $P \in \indextab{\lambda}$. 
\item Let $\psi_P$ be an isomorphism from $S^\lambda$ to $S^P$ for 
$P \in \indextab{\lambda}$. 
\item Let $\Delta_Q^P$ be the image $\psi_P(\Delta_Q)$ of the Specht polynomial $\Delta_Q$.
\item Let $\stdtab{\lambda}$ be the set of standard tableaux on $\lambda$.
\item Let $\stdtabp{\lambda}$ be a maximal subset of  $\stdtab{\lambda}$
such that 
$\Set{ \YYoungsub{\mu} \Delta_{Q} | Q \in \stdtabp{\lambda} }$
is $\CC$-linearly independent. \label{youngsubstep}
\item Let  $\stdtabpp{\lambda}$ be a maximal subset of  $\stdtabp{\lambda}$
such that 
$\Set{ \ZZsub{\mu}{k}{l} \Delta_{Q} | Q \in \stdtabpp{\lambda} }$
is $\CC$-linearly independent.
\item Let $\stdtabpp{\lambda} = \Set{Q_1,Q_2,\ldots}$.
Let $\indextab{\lambda} = \Set{P_1,P_2,\ldots, P_d}$. 
Let $\indextabpair{\lambda} = \Set{(P_i,Q_i)| i=1, \ldots , d}$.
\item Let $f$ be 
$\sum_{\lambda \partitionof m}\sum_{(P,Q)\in \indextabpair{\lambda}} \Delta^{P}_{Q}$.
\end{enumerate}
Hence 
\begin{align*}
\varphi_{\ZZsub{\mu}{k}{l}, \ZZsub{\mu}{k}{l} f}&:
\CC[S_m]\ZZsub{\mu}{k}{l}
\ni \ZZsub{\mu}{k}{l} \mapsto \ZZsub{\mu}{k}{l} f \in
R_{\mu}(k;l)\\
\intertext{and} 
\varphi_{\varepsilon\otimes e , \ZZsub{\mu}{k}{l} f}&:
\induce{H_\mu(l)}{S_m}{Z_\mu(k;l)}
\ni \varepsilon\otimes e \mapsto \ZZsub{\mu}{k}{l} f \in
R_{\mu}(k;l)
\end{align*}
 are isomorphisms of $S_m$-modules. 
\end{theorem}

\begin{remark}
In Theorem \ref{procelemma},
we may replace Specht polynomials by some other bases.
We may also skip Step \ref{youngsubstep}.
\end{remark}

\begin{remark}
Any generic element $\ZZsub{\mu}{k}{l} f$ is obtained by the procedure in
Theorem \ref{procelemma}.
\end{remark}

\begin{example}
We consider the case where $\mu=(1,1,1)$.
In this case, $R_{\mu}$ is the coinvariant algebra of $S_3$.
We can decompose $R_\mu$ as follows:
\begin{align*}
R_\mu &= R_0 \oplus R_1 \oplus R_2 \oplus R_3,\\
R_\mu^0&= \CC \Set{\Delta^{\tiny \young(123)}_{\tiny \young(123)}}, \\
R_\mu^1&= \CC \Set{\Delta^{\tiny \young(12,3)}_{\tiny \young(12,3)}, \Delta^{\tiny \young(12,3)}_{\tiny \young(13,2)}}, \\
R_\mu^2&= \CC \Set{\Delta^{\tiny \young(13,2)}_{\tiny \young(12,3)}, \Delta^{\tiny \young(13,2)}_{\tiny \young(13,2)}}, \\
R_\mu^3&= \CC \Set{\Delta^{\tiny \young(1,2,3)}_{\tiny \young(1,2,3)}}, 
\end{align*}
where $\Delta^P_{Q}$ is the image of the higher Specht polynomial
introduced in \cite{aty}.
Hence we have the decompositions
\begin{align*}
R_\mu (0;3)&= \CC \Set{\Delta^{\tiny \young(123)}_{\tiny \young(123)}} \oplus  \CC \Set{\Delta^{\tiny \young(1,2,3)}_{\tiny \young(1,2,3)}}, \\
R_\mu (1;3)&= \CC \Set{\Delta^{\tiny \young(12,3)}_{\tiny \young(12,3)}, \Delta^{\tiny \young(12,3)}_{\tiny \young(13,2)}}, \\
R_\mu (2;3)&= \CC \Set{\Delta^{\tiny \young(13,2)}_{\tiny \young(12,3)}, \Delta^{\tiny \young(13,2)}_{\tiny \young(13,2)}}.
\end{align*}
For a standard tableau $P \in \stdtab{\lambda}$,
the map from the Specht module $S^\lambda$ 
to $\CC \Set{ \Delta^{P}_{Q} | Q \in \stdtab{\lambda}}$ that 
maps the Specht polynomial $\Delta_{Q}$
to the higher Specht polynomial $\Delta^{P}_{Q}$
for $Q \in \stdtab{\lambda}$ is an isomorphism.
In this case, $ \ZZsub{\mu}{k}{3} = \frac{1}{3}(\varepsilon + \omega^k (123) - \omega^{2k}(132))$, 
where $\omega$ is $\zeta_3$.
By direct calculations,
we have
\begin{align*}
\ZZsub{\mu}{0}{3} \Delta_{\tiny \young(123)} &= \Delta_{\tiny \young(123)}, \\
\ZZsub{\mu}{0}{3} \Delta_{\tiny \young(1,2,3)} &= \Delta_{\tiny \young(1,2,3)}, \\
\ZZsub{\mu}{1}{3} \Delta_{\tiny \young(12,3)} &= \frac{1}{3} \left(
  (2 + \omega)\Delta_{\tiny \young(12,3)}    
- (1 + 2\omega)\Delta_{\tiny \young(13,2)}
\right), \\
\ZZsub{\mu}{2}{3} \Delta_{\tiny \young(13,2)} &= \frac{1}{3} \left(
- (1 + 2\omega)\Delta_{\tiny \young(12,3)}
+ (2 + \omega)\Delta_{\tiny \young(13,2)}    
\right).
\end{align*}
Let 
\begin{align*}
f^{(0)}&=\Delta^{\tiny \young(123)}_{\tiny \young(123)} + \Delta^{\tiny \young(1,2,3)}_{\tiny \young(1,2,3)}, \\
f^{(1)}&= \Delta^{\tiny \young(12,3)}_{\tiny \young(12,3)}, \\
f^{(2)}&= \Delta^{\tiny \young(13,2)}_{\tiny \young(13,2)}.
\end{align*}
The polynomials $\ZZsub{\mu}{k}{3} f^{(k)}$ are generic in 
$R_\mu (k;3)$.
Hence 
\begin{align*}
\varphi_{\ZZsub{\mu}{k}{3}, \ZZsub{\mu}{k}{3} f^{(k)}}&:
\CC[S_3]\ZZsub{\mu}{k}{3}
\ni \ZZsub{\mu}{k}{3} \mapsto \ZZsub{\mu}{k}{3} f^{(k)} \in
R_{\mu}(k;3)\\
\intertext{and} 
\varphi_{\varepsilon\otimes e , \ZZsub{\mu}{k}{3} f^{(k)} }&:
\induce{H_{\mu}(3)}{S_m}{Z_\mu(k;3)}
\ni \varepsilon\otimes e \mapsto \ZZsub{\mu}{k}{3} f^{(k)} \in
R_{\mu}(k;3)
\end{align*}
 are isomorphisms of $S_3$-modules. 
\end{example}

\subsection{The Case of Two Rows }\label{tworowsec}
We have the procedure to construct generic elements (Theorem \ref{procelemma}).
In this section, we show Lemma \ref{onlysslemma},
which makes  the procedure simpler.
Then we show Theorem \ref{tworowscase},
which 
explicitly constructs isomorphisms 
in the case of the Springer modules corresponding to  $\mu=(n,n)$.

\begin{definition}
For $\mu\partitionof m$,
we define $\nu_\mu$ that maps $i\in\Set{1,\ldots,m}$
to $r$ if $i$ lies in the $r$-th row of $t_\mu$.
\end{definition}

\begin{lemma}\label{samelemma}
Let $\lambda$ and $\mu$ be  partitions of $m$.
If standard tableaux $Q$ and $Q'$ on $\lambda$ satisfy $\nu_\mu \circ Q = \nu_\mu \circ Q'$, 
then $\YYoungsub{\mu} \Delta_{Q} = \YYoungsub{\mu} \Delta_{Q}$.
\end{lemma}
\begin{proof}
Since standard tableaux $Q$ and $Q'$ on $\lambda$ satisfy $\nu_\mu \circ Q = \nu_\mu \circ Q'$,
there exists $\sigma \in S_\mu$ such that 
$ Q = \sigma Q' $. Hence we have 
$ S_\mu \Set{Q} = S_\mu \Set{Q'}$,
which implies $\YYoungsub{\mu} \Delta_{Q} = \YYoungsub{\mu} \Delta_{Q'}$.
\end{proof}

\begin{lemma}\label{semisstadlemma}
Let $x$ and $y$ be numbers lying in the same row of $t_\mu$.
If $x$ and $y$ lie in the same column of $Q$,
then $(x,y) \in S_\mu$ acts on $\Delta_Q$ by
\begin{gather*}
(x,y) \Delta_Q = -\Delta_Q.
\end{gather*} 
This implies
$\YYoungsub{\mu} \Delta_Q = 0.$
\end{lemma}
\begin{proof}
It follows from a direct calculation that
\begin{gather*}
(x,y) \Delta_Q = -\Delta_Q.
\end{gather*} 
Let $S_\mu= S'_\mu\Set{\varepsilon} \disjointunion S'_\mu \Set{(x,y)}$.
Since
\begin{align*}
\sum_{\sigma \in S_\mu } \sigma\Delta_Q 
&= \sum_{\sigma \in S'_\mu} \sigma\Delta_Q + \sum_{\sigma \in S'_\mu} \sigma(x,y) \Delta_Q   \\
&= \sum_{\sigma \in S'_\mu} \sigma(\varepsilon + (x,y))\Delta_Q   \\
&= \sum_{\sigma \in S'_\mu} \sigma(\Delta_Q - \Delta_Q)  \\
&=0,
\end{align*}
we have
$\YYoungsub{\mu} \Delta_Q = 0.$
\end{proof}

\begin{lemma}\label{onlysslemma}
Let $\lambda$ and $\mu$ be  partitions of $m$.
Let $\stdtab{\lambda}'$ be a maximal subset of  $\stdtab{\lambda}$
such that 
$\Set{ \YYoungsub{\mu} \Delta_{Q} | Q \in \stdtab{\lambda}' }$
is $\CC$-linearly independent.
Let  $Q, Q'\in \stdtab{\lambda}'$. 
Then $\nu_\mu \circ Q$ is a semi-standard tableau on $\lambda$ with weight $\mu$.
If $Q\neq Q'$, then $\nu_\mu \circ Q \neq \nu_\mu \circ Q'$. 
\end{lemma}
\begin{proof}
By Lemma \ref{semisstadlemma},
 $\nu_\mu  \circ Q$ is a semi-standard tableau on $\lambda$.
It is clear that the weight of  $\nu_\mu \circ Q$  is $\mu$.
Lemma \ref{samelemma} implies
 $\nu_\mu  \circ Q \neq \nu_\mu  \circ Q$ if $Q\neq Q'$.
\end{proof}

By Lemma \ref{onlysslemma}, 
it is enough for  Step \ref{youngsubstep} of the procedure of Theorem 
\ref{procelemma} 
that we consider only the standardizations of 
semi-standard tableaux with weight $\mu$.

Lastly we consider the case where $m=2n$ and $\mu=(n,n)$.

\begin{theorem}\label{tworowscase}
Let $t'_{(2n-k,k)}$ be the standard tableau on $(2n-k,k)$
whose entries in the second row are $\Set{n+1, n+2, \ldots , n+k}$.
Let
\begin{align*}
f^{(0)}&=\YYoungsub{(n,n)}  \sum_{k=0,2,4,6,\ldots}^{k \leq n}
 \Delta_{t'_{(2n-k,k)}}, 
\intertext{and let}
f^{(1)}&=\YYoungsub{(n,n)}  \sum_{k=1,3,5,7,\ldots}^{k \leq n}
 \Delta_{t'_{(2n-k,k)}} 
\end{align*}
where $\Delta_P^Q$ are the images of Specht polynomials.

Then the polynomial $F^{(i)}=\ZZsub{(n,n)}{i}{2} f^{(i)}$ is  generic
 in $R_{(n,n)}(i;2)$. Hence 
$\varphi_{\varepsilon\otimes e, F^{(i)}}: \induce{H_{(n,n)}(2)}{S_{2n}}{Z_{n,n}(i;2)} \to R_{(n,n)}(i;2)$ %,
and  $\varphi_{\ZZsub{(n,n)}{i}{2}, F^{(i)}}:\CC[S_{2n}]\ZZsub{(n,n)}{i}{2}\to R_{(n,n)}(i;2)$
are isomorphisms. 
\end{theorem}
\begin{proof}
For all $k$,
tableaux that are standardizations of 
semi-standard tableaux with weight $(n,n)$ on $(2n-k,k)$
are only $t'_{(2n-k,k)}$.
Since $a_{(n,n),2}=(1,1+n)(2,2+n)\cdots(n,2n)$,
it follows from a direct calculation that $a_{(n,n),2} \Delta_{t'_{(2n-k,k)}} = (-1)^k \Delta_{t'_{(2n-k,k)}}$.
Hence we have the theorem.
\end{proof}

\end{document}